\documentclass[usenatbib]{mn2e}

\usepackage[english]{babel}
\usepackage{hyperref}
\usepackage{graphicx}
\usepackage{enumerate}
\usepackage{amsmath, amsbsy, amsfonts}
\usepackage{multicol}
 
\def\Zbb{\mathbb{Z}}
\def\Rbb{\mathbb{R}}
\def\Tbb{\mathbb{T}}
\def\Gscr{\mathcal{G}}
\def\Rscr{\mathcal{R}}
\def\Vscr{\mathcal{V}}
\def\phi{\varphi}
\def\rho{\varrho}
\def\eps{\varepsilon}

\title{Effective resonant stability of Mercury}

\author[M. Sansottera, C. Lhotka and A. Lema{\^\i}tre]%
       {M. Sansottera\thanks{Email: marco.sansottera@unimi.it},%
        C. Lhotka\thanks{Email: christoph.lhotka@oeaw.ac.at},%
        A. Lema{\^\i}tre\thanks{Email: anne.lemaitre@unamur.be},%
        \\
        $*$Dipartimento di Matematica, Universit\`a degli Studi di Milano,
        Via Cesare Saldini, 50, 20133 Milano, Italy,
        \\
        $\dag$Space Research Institute, Austrian Academy of Sciences, 
        Schmiedlstrasse, 6, 8042, Graz, Austria,
        \\
        $\ddag$ naXys, Universit\'e de Namur,
        Rue de Bruxelles, 61, 5000 Namur, Belgium.}

\begin{document}

\volume{452}
\pagerange{4145--4152}
\pubyear{2015}

\maketitle

\label{firstpage}

\begin{abstract}
Mercury is the unique known planet that is situated in a 3:2
spin-orbit resonance nowadays. Observations and models converge to the
same conclusion: the planet is presently deeply trapped in the
resonance and situated at the Cassini state $1$, or very close to it.
We investigate the complete non-linear stability of this equilibrium,
with respect to several physical parameters, in the framework of
Birkhoff normal form and Nekhoroshev stability theory. We use the same
approach adopted for the 1:1 spin-orbit case, published in
\cite{SanLhoLem-2014}, with a peculiar attention to the role of
Mercury's non negligible eccentricity.  The selected parameters are
the polar moment of inertia, the Mercury's inclination and
eccentricity and the precession rates of the perihelion and node.  Our
study produces a bound to both the latitudinal and longitudinal
librations (of 0.1 radians) for a long but finite time (greatly
exceeding the age of the solar system).  This is the so-called
effective stability time.  Our conclusion is that Mercury, placed
inside the 3:2 spin-orbit resonance, occupies a very stable position
in the space of these physical parameters, but not the most stable
possible one.
\end{abstract}


\section{Introduction}
Mercury is the target of the BepiColombo mission, one of ESA's
cornerstone space missions, carried out in collaboration with the
Japanese Aerospace Agency (JAXA).  The spacecraft will be launched in
2017 and the orbit phase around Mercury is planned in 2024 (please
refer to the BepiColombo webpage on the ESA website,
\url{http://sci.esa.in/bepicolombo}, for updated information).
Mercury has a peculiar feature: it is the only planet in the Solar
System that is locked in a spin-orbit resonance, and the only object
in the Solar System trapped in a 3:2 resonance.  Indeed, the Moon and
most of the regular satellites of the giant planets are found in 1:1
spin-orbit resonance.

The unique situation of Mercury can partly be explained by its large
orbital eccentricity, see, e.g. \cite{ColSha-1966}, \cite{CorLas-2004}
and \cite{CelLho-2014}.  A more realistic tidal model has been used in
\cite{Noyelles201426}, where the authors demonstrate that capture in
3:2 resonance is possible on much shorter time-scales than previously
thought.

These spin-orbit locked positions are usually named the (generalized)
Cassini states ($1$ to $4$) and can be suitably described in terms of
celestial dynamics, see, e.g. \cite{Colombo-1966}, \cite{Peale-1969},
\cite{Beletskii-1972} and \cite{Ward-1975}. For an extension of
  the theory including the polar motion, see also
  \cite{BouKinSou-2003}.

The presence of a spin-orbit resonance allows to link the
observational data of the orbital and rotational states with interior
structure models, this is, for Mercury, the so-called Peale's
experiment, see \cite{Peale-1976}.  

The first observational confirmation of the Cassini State 1 for
Mercury is due to \cite{Margot04052007}. Earth-based radar
observations have confirmed its presence with high accuracy, see,
e.g., \cite{MarEtal-2012}, where the authors demonstrate that the
angle between the spin axis and the orbit normal, commonly referred to
as obliquity, is consistent with the equilibrium hypothesis.

The same value of the obliquity was also determined by observational
results from the NASA space mission MESSENGER that has also validated
the 3:2 spin-orbit resonance to high accuracy, see
\cite{MazEtal-2015}. The most recent value of the obliquity is
$\varepsilon=2.06\pm0.16$ arcmin, the orbital period is given by
$T_o=87.969216879$ days $\pm 6$ seconds and the spin period is
$T_s=58.64616\pm0.000011$ days, that gives a ratio of about 3:2 to
great accuracy.

It has been shown in \cite{Peale-2005} that small free
oscillations around the exact equilibrium of the spin-orbit resonance
are damped due to dissipative forces (mainly tidal effects and
core-mantle friction) on a timescale of $10^{5}$ years. Any effect,
that may bring Mercury away from exact spin-orbit resonance, like
impacts, or any mechanism that may change Mercury's internal mass and
momenta distribution, will be counteracted on relatively short time
scales. It is therefore reasonable to think that Mercury is currently
situated at exact resonance or very close to it.

The question arises concerning the stability of the equilibrium itself
in absence of dissipative forces to separate the influence of
conservative non-linear effects from dissipative ones that act on shorter
time-scales.  The stability of the spin-orbit
resonances has been numerically investigated in, e.g.,
\cite{CelChi-2000}, \cite{CelVoy-2010} and \cite{Lhotka-2013} by means
of stability maps.  The non-linear stability of the Cassini states, in the 1:1
spin-orbit resonance, has been investigated in detail in
\cite{SanLhoLem-2014}, by means of normal forms and Nekhoroshev type
estimates.  However, such an investigation of the non-linear stability,
is still lacking for Mercury.  We aim to carry out this study in the
present paper.  In particular, we are interested in the long-term
non-linear stability: our goal is to produce a bound to both the
latitudinal and longitudinal librations over long but finite times,
namely an \emph{effective stability time}.

Let us stress that in the present work we consider a realistic model
in the {\it mathematical} sense.  Indeed we are able to obtain
significative analytic estimates on the stability time using the {\it
  real} physical parameters.  However, in order to obtain a better {\it
  physically} realistic model, one should also take into account
dissipative effects and planetary perturbations.

The long-term stability of perturbed proper rotations (rotations about
a principal axis of inertia) in the sense of Nekhoroshev has 
been shown in \cite{BenettinFassoGuzzo2004}.  Furthermore, the authors
suggest that their results may possibly be extended to the case of
spin-orbit resonances.  With our study we are able to demonstrate the
long-term stability for motions that are trapped in a spin-orbit
resonance and in particular the possible application of Nekhoroshev
theory to Mercury.

We are able to give a definitive answer: the generalized Cassini state
$1$, realized in terms of the 3:2 resonance, is practically stable on
long time scales. However, we also demonstrate that the actual
position of Mercury, in the parameter space, is placed in a very
stable region, while not the most stable one.  Indeed, altering
some physical parameters, namely the polar moment of inertia factor,
the inclination, the eccentricity and the precession rates of the
perihelion and node, we found that the stability may change from a
marginal amount to orders of magnitude, depending on the quantity of
interest.

The present paper represents an extension of our previous work,
developed for the 1:1 spin-orbit problems, see \cite{SanLhoLem-2014},
to the case of the 3:2 resonance and in particular to Mercury. The
presence of a non negligible value for the eccentricity is a new
aspect of the model to take into account.  The mathematical basis of
our work is represented by the Birkhoff normal
form (\citeyear{Birkhoff-1927}) and the Nekhoroshev
theory (\citeyear{Nekhoroshev-1977, Nekhoroshev-1979}).

Our approach is reminiscent of similar works on the same line,
see, e.g., \cite{GioLocSan-2009, GioLocSan-2010} and \cite{SanLocGio-2011}, in
which the authors gave an estimate of the long-time stability for the
Sun-Jupiter-Saturn system and the planar Sun-Jupiter-Saturn-Uranus
system, respectively.

This work is organized as follows: we introduce the spin-orbit model
and its Hamiltonian formulation in Section~\ref{sec:1}.  In
Section~\ref{sec:2} we describe an algorithm for the evaluation
of the effective stability time via Birkhoff normal form.  The
application of our study to Mercury is presented in
Section~\ref{sec:3}, while the physical interpretations and the
possible extensions are reported in Section~\ref{sec:4}.

\section{The model}\label{sec:1}
We consider Mercury as a triaxial rigid body whose principal moments
of inertia are $A$, $B$ and $C$, with $A\leq B<C$.  We denote by $m$
and $R_{\rm e}$, the mass and equatorial radius of Mercury,
respectively; and by $m_0$ the mass of the Sun.

We closely follow the notation adopted in \cite{SanLhoLem-2014}, that
was also used in some previous studies on the same subject, see,
e.g., \cite{HenSch-2004} for a general treatment of synchronous
satellites, \cite{DHoLem-2004} and \cite{LemDHoRam-2006a} for Mercury,
\cite{NoyLemVie-2008} for the study of Titan, \cite{Lhotka-2013}
for a symplectic mapping model and \cite{NoyLho-2013} for an
investigation concerning the obliquity of Mercury during the
BepiColombo space mission.  Thus we refer to the quoted works for a
detailed exposition and we just report here the key points so as to
make the paper quite self-contained.  First, we briefly recall how to
express the Hamiltonian in the Andoyer-Delaunay set of coordinates,
then we introduce the simplified spin-orbit model that represents the
basis of our study.

\subsection{Reference frames}\label{sbs:ReferenceFrames}
The usual description of the spin-orbit motion requires four basic
reference frames having their common origin in the center of mass of
Mercury, see \cite{DHoLem-2004}.  Namely, we need: (i)~the
\emph{inertial frame}, $(X_0, Y_0, Z_0)$, with $X_0$ and $Y_0$ lying
in the ecliptic (or Laplace) plane; (ii)~the \emph{orbital frame}, $(X_1, Y_1,
Z_1)$, with $Z_1$ perpendicular to the orbit plane; (iii)~the
\emph{spin frame}, $(X_2, Y_2, Z_2)$, with $Z_2$ pointing to the
spin axis direction and $X_2$ to the ascending node of the equatorial
plane in the ecliptic plane; (iv)~the \emph{body frame}, (or
\emph{figure frame}), with $Z_3$ pointing into the direction of the
axis of greatest inertia and $X_3$ of smallest inertia.  In order to
link the different frames we define $\nu_{01}$ as the direction along
the ascending node of $(X_0,Y_0)$, in the plane $(X_1,Y_1)$, and
$\nu_{23}$ as the direction of the ascending node between the planes
$(X_2,Y_2)$ and $(X_3,Y_3)$.

We introduce two Euler angles: (i)~the \emph{inertial obliquity}, $K$,
that is the angle between the axes $Z_0$ and $Z_2$; (ii)~the
\emph{wobble}, $J$, between the axes $Z_2$ and $Z_3$.  Moreover, we
define the angles for the rotational motion: (i)~the \emph{spin node},
$h_s$, between $X_0$ and $X_2$, measured in the plane $(X_0,Y_0)$;
(ii)~the \emph{figure node}, $g_s$, between $\nu_{23}$ and $X_3$,
measured in the plane $(X_2,Y_2)$; (iii)~the \emph{rotation angle},
$l_s$, between $\nu_{23}$ and $X_3$, measured in the plane $(X_3,
Y_3)$.

For the orbital dynamics, we introduce: (i)~the \emph{longitude of the
  ascending node}, $\Omega$, that gives the direction of $\nu_{01}$
measured in the plane $(X_0,Y_0)$; (ii)~the \emph{inclination}, $i$,
being the angle between the axes $Z_0$ and $Z_1$; (iii)~the
\emph{perihelion argument}, $\omega$, that defines the direction of
the pericenter $X_1$ in the plane $(X_1,Y_1)$.

\begin{figure}
\centering
\vskip5pt
\includegraphics[width=.35\linewidth]{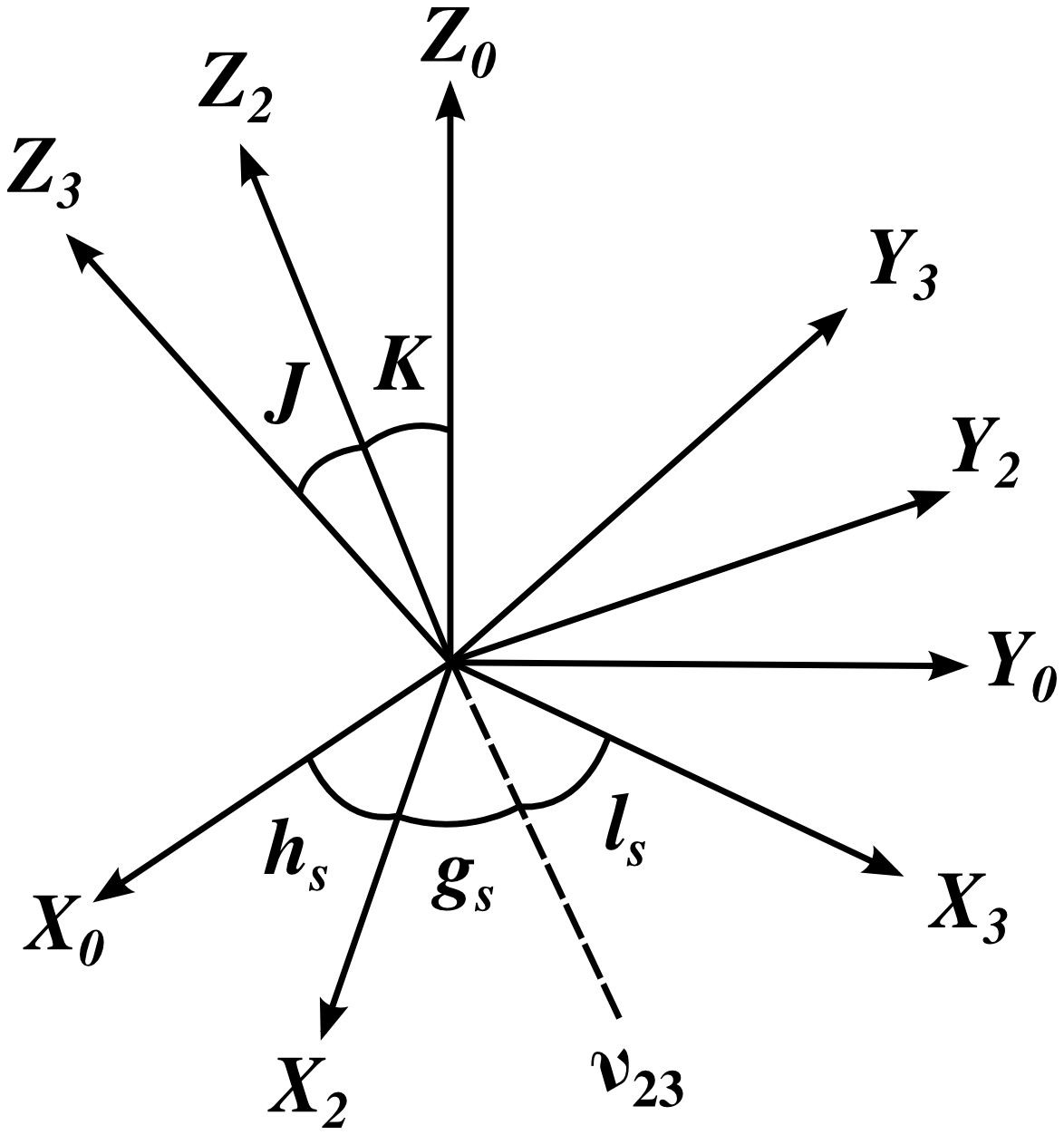}
\hspace{30pt}
\includegraphics[width=.35\linewidth]{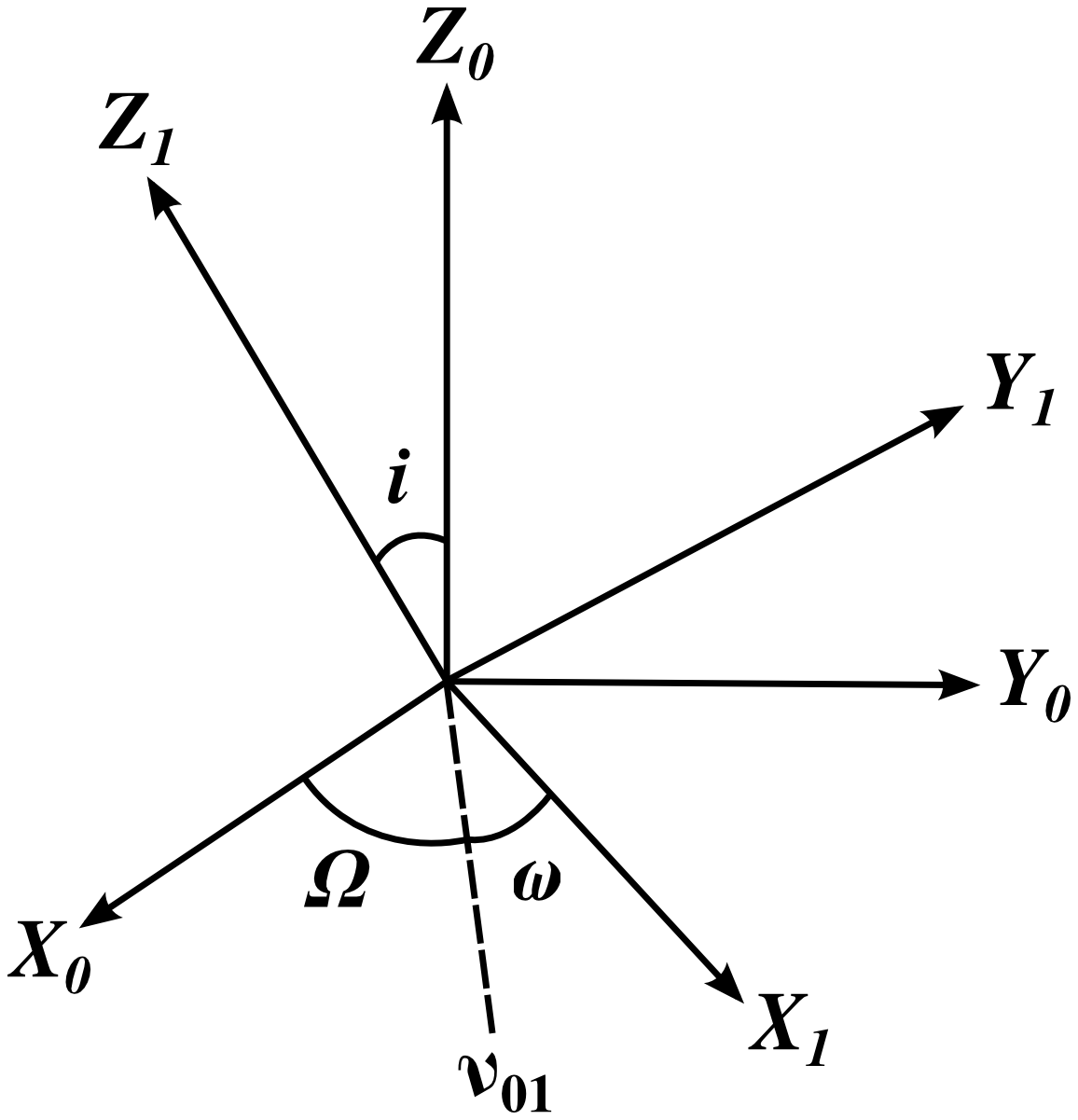}
\caption{The four reference frames and the relevant angles related to
  the Andoyer (left) and Delaunay (right) canonical variables.  See
  the text for more details.}\label{fig:sop}
\end{figure}

We report in Figure~\ref{fig:sop} the four reference frames and the
angles defined above. On its basis we introduce the Andoyer-Delaunay 
canonical variables.

\subsection{Andoyer-Delaunay variables}
We use the \emph{modified Andoyer} variables in order to describe the
rotational motion:
\begin{equation}
\begin{aligned}
L_1 &= G_s     \ , &l_1 &=  l_s+g_s+h_s\ ,\\
L_2 &= G_s-L_s \ , &l_2 &= -l_s\ ,\\
L_3 &= G_s-H_s \ , &l_3 &= -h_s\ ,\\
\end{aligned}
\label{eq:Andoyer}
\end{equation}
where $G_s$ denotes the norm of the spin angular momentum,
$L_s=G_s\cos{J}$ and $H_s=G_s\cos{K}$.

For the orbital motion, we use the \emph{modified Delaunay} variables,
namely
\begin{equation}
\begin{aligned}
L_4 &= L_o     \ , &l_4 &=  M+\omega+\Omega\ ,\\
L_5 &= L_o-G_o \ , &l_5 &= -\omega-\Omega\ ,\\
L_6 &= G_o-H_o \ , &l_6 &= -\Omega\ ,
\end{aligned}
\label{eq:Delaunay}
\end{equation}
with $L_o = m\sqrt{\mu a}$, $G_o = L_o\sqrt{1-e^2}$, and $H_o =
G_o\cos{i}$.  As usual, $\mu=\Gscr(m+m_0)$ where $\Gscr$ is the
gravitational constant, $a$ is the semi-major axis, $e$ the
eccentricity, $i$ the inclination and $M$ the mean anomaly.  Let us
remark that $L_5$ and $L_6$ are related to the eccentricity, $e$, and
inclination, $i$, respectively. From now on, we assume that the
spin-axis is aligned with the line of figure, i.e., $J=0$, and that we
can neglect motion related to the conjugated variables $(L_2,l_2)$.
Indeed, the polar motion of Mercury is assumed to be very
  small, see, e.g., \cite{NoyDufLem-2010}. Taking into account the
3:2 spin-orbit resonance, we introduce the resonant variables
$\Sigma=(\Sigma_1,\Sigma_3)$, $\sigma=(\sigma_1,\sigma_3)$ by using
the generating function
$$
S_{3:2}=\Sigma _1\left(l_1-\frac{3}{2}\big(l_4+l_5\big)+l_5\right)+
\Sigma _3\left(l_3-l_6\right) \ .
$$
Thus we get
\begin{align}
\Sigma _1&=L_1 \ , & \sigma _1&=l_1-\frac{3}{2}\big(l_4+l_5\big)+l_5 \ ,
\nonumber\\ 
\Sigma _3&=L_3 \ , &\sigma _3&=l_3-l_6 \ . 
\label{resvar}
\end{align}
Here, $\sigma _1$ and $\sigma _3$ refer to the longitudinal and
latitudinal librations, respectively, around the exact resonant state,
placed at $\sigma _1 = \sigma _3 = 0$.

In our model, the orbital ellipse of Mercury is assumed to be frozen
but uniformly precessing due to the interactions with the other planets
of the solar system.  We keep fixed the semi-major axis, the eccentricity
and the inclination: this corresponds to fixing the values of $L_4$,
$L_5$ and $L_6$.  Denoting by $n$ the mean motion, by $\dot\omega$ the
mean precession rate of the argument of perihelion, by $\dot\Omega$
the mean regression rate of the ascending node and, without loss of
generality, setting $t_0=0$ as the time of the perihelion passage of
Mercury, we have the trivial equations of motion for the conjugated
angles: $l_4=n t$, $l_5=\dot\omega t$ and $l_6=\dot\Omega t$.
In this setting, the generating function $S_{3:2}$ is a
  function of time.  Therefore, in order to express the Hamiltonian in
  the resonant variables, we have to add also the time derivative of
  the generating function, $\partial S_{3:2}/\partial t$.

\subsection{Hamiltonian formulation}
Let us denote by ${\mathcal H}_s$ and ${\mathcal H}_o$ the rotational and
orbital kinetic energy, respectively, then the kinetic part reads

\begin{equation}
\begin{split}
{\mathcal T }&=\mathcal{H}_s+\mathcal{H}_o+
\frac{\partial S_{3:2}}{\partial t}   \nonumber \\
&=\frac{\Sigma_1^2}{2C}
-\frac{3}{2}n\Sigma_1
-\Sigma_1\dot{\omega}
+\left(\Sigma_3-\Sigma_1\right)\dot{\Omega}   \ ,
\end{split}
\end{equation}
where $C$ is the Mercury's largest moment of inertia.

The dominant contribution of the gravitational potential is
mainly due to two spherical harmonics: $C_{20}=\frac{A+B-2C}{2 m R_{\rm e}^2}$ and
$C_{22}=\frac{B-A}{4  m R_{\rm e}^2}$.  Thus we consider the
potential energy
\begin{equation}
\label{potful}
{\mathcal V} = -\frac{\Gscr m m_0 R_{\rm e}^2}{r^3} \Bigl(
C_{20}P_{20} (\sin\phi) +
C_{22}P_{22} (\sin\phi)\cos 2 \lambda
\Bigr) \ ,
\end{equation}
where $P_{20}$, $P_{22}$ are the Legendre polynomials and $\phi$,
$\lambda$ are the co-latitude and longitude of Mercury, respectively,
with $\lambda$ measured east-wards from $X_2$ ($X_3$).  The expression
of the potential ${\mathcal V}$, defined in the figure frame, into the
inertial frame is straightforward.  Following the approach developed
in \cite{NoyLho-2013}, we write ${\mathcal V}$ in terms of the
resonant variables and perform an average over the mean orbital
longitude $l_4$.  We denote by $\langle {\mathcal V}\rangle$ the
averaged potential and we refer to the previous paper for all the
details.  We just stress that here, as an extension to
\cite{NoyLho-2013}, we consider an expansion in the eccentricity up to
order 8, while we only keep the time-independent harmonics in the
definition of ${\langle \Vscr\rangle}$, namely the ${\langle
  \Vscr_{20}\rangle}$ and ${\langle \Vscr_{22}\rangle}$ terms, defined
there.

Denoting by $\langle T \rangle$ the averaged kinetic energy and
introducing the averaged Hamiltonian ${\langle\mathcal
  H\rangle=\langle T \rangle+\langle \Vscr\rangle}$, the equations of
motion read
\begin{equation}
\label{canequ}
\dot \Sigma = - \frac{\partial {\langle\mathcal H\rangle}}{\partial\sigma} \ , \quad
\dot \sigma = + \frac{\partial {\langle\mathcal H\rangle}}{\partial\Sigma} \ .
\end{equation}
Setting $\sigma=0$ we look for the equilibrium
$\Sigma^*=\left(\Sigma_1^*,\Sigma_3^*\right)$ by solving the equations
\begin{equation}
\label{equilib}
f_1\left(\Sigma\right)\equiv
\frac{\partial {\langle\mathcal H\rangle}}{\partial \Sigma_1}\bigg|_{\sigma=0}\hskip-5pt=0 \ , \quad
f_2\left(\Sigma\right)\equiv
\frac{\partial {\langle\mathcal H\rangle}}{\partial \Sigma_2}\bigg|_{\sigma=0}\hskip-5pt=0 \ ,
\end{equation}
namely the Cassini state 1.

Let us remark that, expressing $f_1\left(\Sigma\right)$ and
$f_2\left(\Sigma\right)$ as functions of $(G_s,K)$ one obtains an
implicit formula for the obliquity $\varepsilon=i-K$ in terms of the
system parameters, see \cite{NoyLho-2013}.  Precisely, setting
$c=C/(mR_{\rm e}^2)$, at the equilibrium, the following equation holds
true
\begin{equation}
c=\frac{
n \sin (\varepsilon ) \left(C_{20} H_{20} \cos (\varepsilon )+C_{22}
   H_{22} (\cos (\varepsilon )+1)\right)
}
{
\dot{\Omega } \sin (i-\varepsilon ) \left(\frac{2 \dot{\Omega } \cos
   (i-\varepsilon )}{3 n}+\frac{2 \dot{\omega }}{3 n}+1\right)
} \ ,
\label{myFor}
\end{equation}
where $H_{20}$ and $H_{22}$, truncated at order 8, are given by
\begin{align}
H_{20}&=
-1-
\frac{3 e^2}{2}-
\frac{15 e^4}{8}-
\frac{35 e^6}{16}-
\frac{315 e^8}{128}
\nonumber \ , \\
H_{22}&=
\frac{7e}{2}-
\frac{123 e^3}{16}+
\frac{489 e^5}{128}-
\frac{1763 e^7}{2048} \ .
\end{align}
This relation will be useful for the physical interpretation of the
results in terms of fixed $\varepsilon$.  Let us remark that equation
\eqref{myFor} represents a generalization of the one in
\cite{Peale1981143}.

\section{Stability at the Cassini state}\label{sec:2}
We now aim to study the stability properties of the Cassini state.  To
be more specific, our goal is to give an estimate of the effective
stability time around the equilibrium point. Hereinafter we follow
the exposition given in \cite{SanLhoLem-2014} for the 1:1 spin-orbit problem.

We perform a translation in order to put the equilibrium at the origin
and an expansion of the averaged Hamiltonian in power series of
$(\Sigma,\sigma)$, namely
\begin{equation}
H(\Sigma,\sigma) = H_0(\Sigma,\sigma) + \sum_{j>0} H_j(\Sigma,\sigma)\ ,
\label{eq:Hsigma}
\end{equation}
where $H_j$ is an homogeneous polynomial of degree $j+2$ in
$(\Sigma,\sigma)\,$.  In the latter equation the quadratic term,
$H_0$, has been put apart from the other terms of the Hamiltonian in
view of its relevance in the perturbative scheme.

\subsection{The untangling transformation}\label{sbs:diagonalization}
The quadratic part of the Hamiltonian reads
\begin{align*}
H_0(\Sigma,\sigma) = &\mu_{{\scriptscriptstyle \Sigma_1 \Sigma_1}}
\Sigma_1^2 +2\mu_{{\scriptscriptstyle \Sigma_1 \Sigma_3}} \Sigma_1
\Sigma_3 +\mu_{{\scriptscriptstyle \Sigma_3 \Sigma_3}}
\Sigma_3^2\\ &+\mu_{{\scriptscriptstyle \sigma_1 \sigma_1}} \sigma_1^2
+2\mu_{{\scriptscriptstyle \sigma_1 \sigma_3}} \sigma_1 \sigma_3
+\mu_{{\scriptscriptstyle \sigma_3 \sigma_3}} \sigma_3^2\ .
\end{align*}
We perform the so-called {\it untangling transformation},
see \cite{HenLem-2005}, that permits to get rid of the mixed terms.
Thus, in the new coordinates, $H_0$ takes the form
\begin{equation*}
H_0(\Sigma',\sigma') =
\mu'_{{\scriptscriptstyle \Sigma_1' \Sigma_1'}} \Sigma_1'^2
+\mu'_{{\scriptscriptstyle \sigma_1' \sigma_1'}} \sigma_1'^2
+\mu'_{{\scriptscriptstyle \Sigma_3' \Sigma_3'}} \Sigma_3'^2
+\mu'_{{\scriptscriptstyle \sigma_3' \sigma_3'}} \sigma_3'^2\ .
\label{eq:untangling}
\end{equation*}
Let us remark that if both $\mu'_{{\scriptscriptstyle \Sigma_1'
    \Sigma_1'}}\mu'_{{\scriptscriptstyle \sigma_1' \sigma_1'}}$ and
$\mu'_{{\scriptscriptstyle \Sigma_3'
    \Sigma_3'}}\mu'_{{\scriptscriptstyle \sigma_3' \sigma_3'}}$ are
positive, as it happens in our case, the quadratic part of the
Hamiltonian represents a couple of harmonic oscillators.

It is now useful to perform a rescaling and introduce the polar
coordinates,
\begin{equation}
\begin{aligned}
\Sigma_1' &= \sqrt{2U_1/U_1^*} \cos(u_1)\ , \qquad\sigma_1' = \sqrt{2 U_1 U_1^*} \sin(u_1)\ ,\\
\Sigma_3' &= \sqrt{2U_3/U_3^*} \cos(u_3)\ , \qquad\sigma_3' = \sqrt{2 U_3 U_3^*} \sin(u_3)\ ,
\end{aligned}
\label{eq:actang}
\end{equation}
where
$$
U_1^* = \sqrt{\mu'_{{\scriptscriptstyle \Sigma_1' \Sigma_1'}}/\mu'_{{\scriptscriptstyle \sigma_1' \sigma_1'}}}\ ,
\qquad\hbox{and}\qquad
U_3^* = \sqrt{\mu'_{{\scriptscriptstyle \Sigma_3' \Sigma_3'}}/\mu'_{{\scriptscriptstyle \sigma_3' \sigma_3'}}}\ .
$$
Thus, the quadratic part of the Hamiltonian is expressed in
action-angle variables as
$$
H_0 = \omega_{u_1} U_1 + \omega_{u_3} U_3\ ,
$$ where $\omega_{u_1}$ and $\omega_{u_3}$ are the frequencies of the
angular variables $u_1$ and $u_3$, respectively.  Again, we use
the shorthand notations $U=(U_1,U_3)\,$, $u=(u_1,u_3)$ and
$\omega_{u}=(\omega_{u_1},\omega_{u_3})\,$.

In these new coordinates, the transformed Hamiltonian can be
expanded in Taylor-Fourier series and reads
\begin{equation}
H^{(0)}(U,u) = \omega_{u}\cdot U
+ \sum_{j>0} H^{(0)}_{j}(U,u)\ ,
\label{eq:H_Uu}
\end{equation}
where the terms $H_j$ are homogeneous polynomials of degree $j/2+1$ in
$U$, whose coefficients are trigonometric polynomials in the angles
$u\,$.

Mathematically speaking, the Hamiltonian~\eqref{eq:H_Uu} describes a 
perturbed system of harmonic oscillators, where the perturbation is 
proportional to the distance from the equilibrium.  We now aim to 
investigate the stability around this equilibrium in 
the light of Nekhoroshev theory, introducing the so-called 
\emph{effective stability time}.

\subsection{Effective stability via Birkhoff normal form}
Following a quite standard approach, we first construct the Birkhoff
normal form for the Hamiltonian~\eqref{eq:H_Uu} and then give an
estimate of the stability time.

The Hamiltonian is in normal form at order $r$ if
\begin{equation}
H^{(r)}(U,u) = Z_0(U)
+\ldots+Z_r(U)+\sum_{s>r} \Rscr^{(r)}_{s}(U,u)\ ,
\label{eq:H_r}
\end{equation}
where $Z_s$, for $s=0,\,\ldots,r\,$, is a homogeneous polynomial of
degree $s/2+1$ in $U$ and in particular  is zero for odd $s$. The
unnormalized remainder terms $\Rscr^{(r)}_s$, for $s>r$, are
homogeneous polynomials of degree $s/2+1$ in $U$, whose coefficients
are trigonometric polynomials in the angles $u\,$.

It is well known, see, e.g., \cite{Giorgilli-1988}, that the Birkhoff
normal form at any finite order $r$ is convergent in some neighborhood
of the origin, but the analyticity radius shrinks to zero when the
order $r\to\infty\,$.  Therefore, we look for stability over a finite
time, possibly long enough with respect to the lifetime of the system.
More precisely, we want to produce quantitative estimates that allow
to give a lower bound of the stability time.

We pick two positive numbers $R_1$ and $R_3$, and consider a polydisk
$\Delta_{\rho R}$ centered at the origin of $\Rbb^2$, defined as
$$
\Delta_{\rho R} = \left\{
U\in \Rbb^2 : |U_j| \leq \rho R_j\ ,\ j=1,3
\right\}\ ,
$$
$\rho>0$ being a parameter. We consider a function
$$
f_s(U,u)=\sum_{|l|=s+2,k\in\Zbb^2} f_{l,k}
\,U^{l/2} {\sin\atop\cos} (k\cdot u)\ ,
$$
which is a homogeneous polynomial of degree $s/2+1$ in
the actions $U$ and depends on the angles $u$.  We define the quantity
$|f_s|_R$ as
$$
|f_s|_R = \sum_{|l|=s+2,k\in\Zbb^2} |f_{l,k}| R_1^{l_1/2} R_3^{l_2/2}\ .
$$
Thus we get the estimate
$$
|f(U,u)| \leq |f|_R\, \rho^{s/2+1}\ ,
\quad\hbox{for }
U\in\Delta_{\rho R}\,,\ u\in\Tbb^2\ .
$$

Given $U(0)\in\Delta_{\rho_0 R}$, with $\rho_0 < \rho$, we have
$U(t)\in\Delta_{\rho R}$ for $t\leq T$, where $T$ is the escape time
from the domain $\Delta_{\rho R}$.  The Hamiltonian~\eqref{eq:H_r} is
in Birkhoff normal form up to order $r$, thus we have
$$
|\dot U| \leq |\{ U, H^{(r)} \}| = \sum_{s>r} |\{ U, \Rscr^{(r)}_s \}|
\leq d |\{ U, \Rscr^{(r)}_{r+1} \}|_R\, \rho^{r/2+1}\ ,
$$
with $d\geq1$.  In fact, after having set $\rho$ smaller than the
convergence radius of the remainder series, $\Rscr^{(r)}_s$ for $s>r$,
the above inequality holds true for some value $d$.

The latter equation allows us to find a lower bound for the escape
time from the domain $\Delta_{\rho R}$, namely the time when
physical librations exceed the given threshold, 
\begin{equation}
\tau(\rho_0, \rho, r) = 
\frac{\rho-\rho_0}{d |\{ U, \Rscr^{(r)}_{r+1} \}|_R\, \rho^{r/2+1}}\ ,
\label{eq:time}
\end{equation}
which, however, depends on $\rho_0$, $\rho$, and $r$.  We stress that
$\rho_0$ is the only physical parameter, being fixed by the initial
data, while $\rho$, and $r$ are left arbitrary.  Indeed, the parameter
$\rho_0$ must be chosen in such a way that the domain $\Delta_{\rho_0
  R}$ contains the initial conditions of the system.  In order to
achieve an estimate of the escape time, $T(\rho_0)$, independent of
$\rho$ and $r$, we introduce
\begin{equation}
T(\rho_0) = \max_{r\geq1}\, \max_{\rho>0}\,\tau(\rho_0, \rho, r)\ ,
\label{eq:timemax}
\end{equation}
which is our best estimate of the escape time.  We define this
quantity as the \emph{effective stability time}.

\begin{figure*}
  \begin{center}
    \resizebox{.8\textwidth}{!}{\input{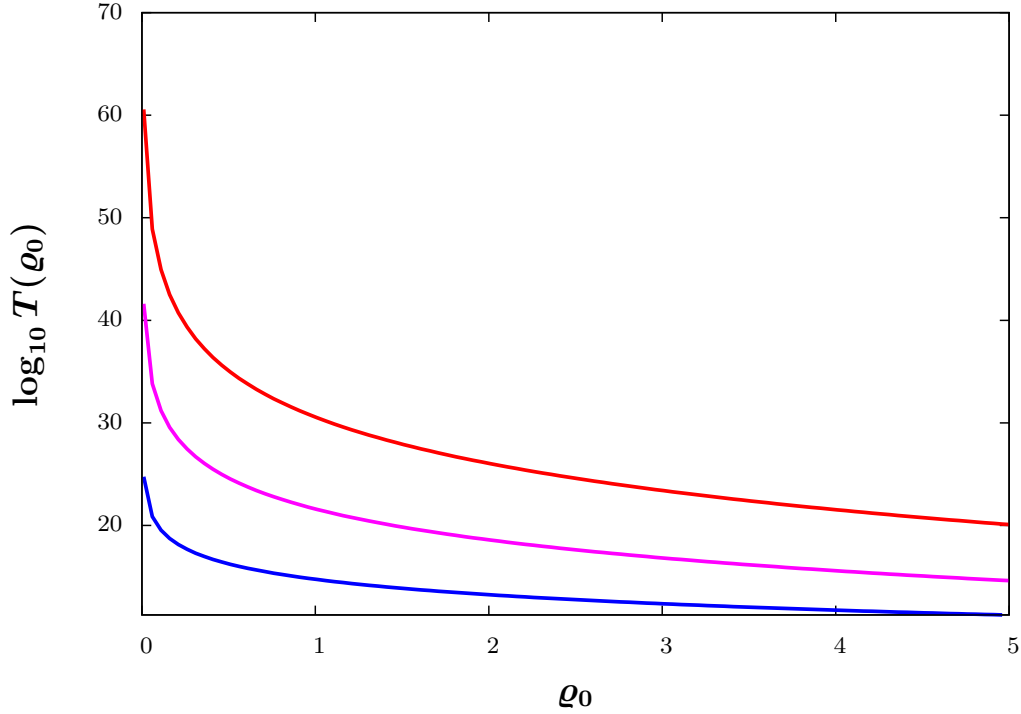}}
    \caption{Estimated effective stability time,
      $T(\rho_0)$, versus normalized distance from the equilibrium
      point, $\rho_0$. The time unit is the year. The three lines
      correspond (from bottom to top) to three different normalization
      orders: r = 10 (blue), r = 20 (pink) and r = 30 (red).}
    \label{fig:time}
  \end{center}
\end{figure*}

\section{Application to Mercury}\label{sec:3}
We now apply the algorithm described in the previous section to
Mercury and evaluate the effective rotational stability time as a
function of some relevant physical parameters.

The expansion of the Hamiltonian function and all the transformations
needed to put the Hamiltonian in the form~\eqref{eq:H_Uu} have been
done using the Wolfram Mathematica software, while the high-order
Birkhoff normal form, up to the order $r=30$, has been computed via a
specific algebraic manipulator, i.e., {{\it X$\rho$\'o$\nu
    o\varsigma$}}, see \cite{GioSan-2012}. 

In the actual calculations we take as reference values the physical
parameters reported in Table~\ref{tab:mercury_param}, where $M$, $m$,
$R$ are taken from \url{http://ssd.jpl.nasa.gov/?planet_phys_par}, we
use $J_2$ and $C_{22}$ given in \cite{Smith13042012}, and $C/m R_{{\rm
    e}}^2$ is obtained by Peale's formula (or \eqref{myFor}) on the basis of the data
provided in \cite{MazEtal-2015}.  Mean orbital elements (J2000) for
$a$, $e$ are taken from
\url{http://ssd.jpl.nasa.gov/txt/p_elem_t2.txt}.  We use the value of
the inclination defined with respect to the Laplace plane,
$i=8.6^\circ$, see, e.g. \cite{Yseboodt2006327}, instead of the
ecliptic.  Finally we compute $\dot\omega$ ($\dot\Omega$) from the
precession period of the perihelion (128 ky), and the regression
period of the ascending node (328 ky), respectively.  Since our
results rely on a qualitative study we removed non-significant digits
from Table~\ref{tab:mercury_param}.

\renewcommand{\arraystretch}{1.1}
\begin{table}
\begin{center}
\caption{Mercury's physical parameters. See the beginning of
  Section~\ref{sec:3} for detailed references for the parameters.}
\begin{tabular}{|c|l|}
\hline
\vphantom{$|^{|^|}$} $M$ & $\phantom{-}1.98843\times10^{+30}\,$ Kg \\
$J_2$ & $\phantom{-}5.031\times10^{-5}$ \\
$C_{22}$ & $\phantom{-}8.088\times10^{-6}$ \\
$C/mR_{\rm e}^2$ & $\phantom{-}3.49\times10^{-1}$ \\
$R_e$ & $\phantom{-}2439.7$ Km\\
$m$ & $\phantom{-}3.30104\times10^{23}$ Kg \\
$a$ & $\phantom{-}5.79091\times10^{+7}\,$ Km\\
$e$ & $\phantom{-}2.05630\times10^{-1}$ \\
$i$ & $\phantom{-}1.50098\times10^{-1}\,$ rad\\
$\dot\omega$ & $\phantom{-}1.34118\times10^{-7}\,$ rad/year\\
$\dot\Omega$ & $-5.23390\times10^{-8}\,$ rad/year\\
$n$ & $\phantom{-}7.1229\times10^{-2}\,$ rad/year\\
\hline
\end{tabular}
\label{tab:mercury_param}
\end{center}
\end{table}

\noindent In Figure~\ref{fig:time} we plot the logarithm of
the effective rotational stability time of Mercury,
$\log_{10}T(\rho_0)$, versus the distance from the equilibrium,
$\rho_0$.  We recall that $\rho_0=0$ corresponds to the Cassini state,
while increasing values of $\rho_0$ allow oscillations around the
equilibrium point.  The highest normalization order, namely $r=30$,
gives the best estimate.  Nevertheless, already at order $r=20$ we
reach an effective stability time greatly exceeding the estimated age
of the Universe, being of the order $10^{10}$, in a domain
$\Delta_{R}$ that roughly corresponds to a libration of $0.1$
radian. Since the actual librations of Mercury around exact resonance are
smaller, we conclude that the spin-orbit coupling of Mercury in a 3:2
spin-orbit resonance is practically stable for the age of the Universe.

\begin{figure}
  \resizebox{\linewidth}{!}{\input{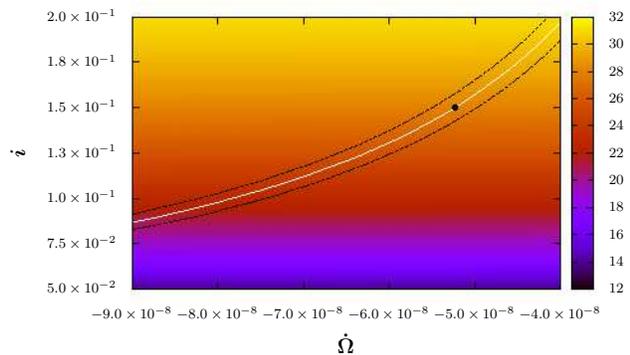}}
  \caption{Stability time: $\dot\Omega$ vs. $i$ (see text for more details).}
  \label{f:om2dot-incl}
\end{figure}

\begin{figure}
  \resizebox{\linewidth}{!}{\input{./time_ecc-incl_U3.tex}}
  \caption{Stability time: $e$ vs. $i$ (see text for more details).}  
\label{f:ecc-incl}
\end{figure}

\begin{figure}
  \resizebox{\linewidth}{!}{\input{./time_crm-incl_U3.tex}}
  \caption{Stability time: $c$ vs. $i$ (see text for more details).}
  \label{f:crm-incl}
\end{figure}

\begin{figure}
  \resizebox{\linewidth}{!}{\input{./time_om1dot-incl_U3.tex}}
  \caption{Stability time: $\dot\omega$ vs. $i$ (see text for more details).}
  \label{f:om1dot-incl}
\end{figure}

\begin{figure}
  \resizebox{\linewidth}{!}{\input{./time_crm-om2dot_U3.tex}}
  \caption{Stability time: $c$ vs. $\dot\Omega$ (see text for more details).}
  \label{f:crm-om2dot}
\end{figure}

\begin{figure}
  \resizebox{\linewidth}{!}{\input{./time_crm-om1dot_U3.tex}}
  \caption{Stability time: $c$ vs. $\dot\omega$ (see text for more details).}
    \label{f:crm-om1dot}
\end{figure}

\begin{figure}
  \resizebox{\linewidth}{!}{\input{./time_om1dot-om2dot_U3.tex}}
  \caption{Stability time: $\dot\Omega$ vs. $\dot\omega$ (see text for more details).}
  \label{f:om1dot-om2dot}
\end{figure}

\begin{figure}
  \resizebox{\linewidth}{!}{\input{./time_crm-ecc_U3.tex}}
  \caption{Stability time: $c$ vs. $e$ (see text for more details).}
  \label{f:crm-ecc}
\end{figure}

\begin{figure}
  \resizebox{\linewidth}{!}{\input{./time_om2dot-ecc_U3.tex}}
  \caption{Stability time: $\dot\Omega$ vs. $e$ (see text for more details).}
  \label{f:om2dot-ecc}
\end{figure}

\begin{figure}
  \resizebox{\linewidth}{!}{\input{./time_om1dot-ecc_U3.tex}}
  \caption{Stability time: $\dot\omega$ vs. $e$ (see text for more details).}
  \label{f:om1dot-ecc}
\end{figure}

\subsection{Sensitivity to physical parameters}\label{sbs:timeparam}
In this section we investigate the dependency of the effective
stability time on the following Mercury's physical parameters: the
polar moment of inertia, $c$, the precession rate of the perihelion
argument, $\dot\omega$, the mean regression rate of the ascending
node, $\dot\Omega$, the eccentricity $e$ and the inclination $i$.
Precisely, we consider $11$ equally spaced different values of each
parameter in the ranges
\halign{\hskip12pt -\ $#$&$\,\in$ $#$\hfil\cr c&
  [0.3,0.4]\ ;\cr \dot\omega & [1\times 10^{-7}, 2\times
    10^{-7}]\ {\rm (rad/day)}\ ;\cr \dot\Omega & [-9\times
    10^{-8},-4\times 10^{-8}]\ {\rm (rad/day)}\ ;\cr e & [0,0.4]\ ;\cr
  i & [0.05,0.2]\ {\rm (rad)}\ .\cr}

The ranges are chosen in order to include possible variations of the
orbital elements due to planetary perturbations, see,
e.g., \cite{Laskar20081}. The choice for polar moment of inertia $c$ is
motiviated to include a variety of possible interior structure models,
from a thin shell to a homogenous sphere.

For each choice of parameters we compute the effective stability time
\eqref{eq:timemax} using the procedure outlined above. As a by-product
we also compute numerically the equilibrium point $\Sigma^*$
using~\eqref{equilib}. From~\eqref{eq:Andoyer} we then get the
corresponding value of the inertial obliquity, namely $K_*$, from
which we obtain the observable obliquity defined as
$\eps_*=K_*-i$. For testing purpose we cross-check $\eps_*$ with
solutions directly obtained from \eqref{myFor}. We notice that the
knowledge of specific values of $\eps$ for different given parameters
allows us to relate the information about the stability time with
observations, i.e. the current observed value for Mercury
$\eps=2.06'$. We present our results by means of contour plots, see
Figures~\ref{f:om2dot-incl}--\ref{f:om1dot-ecc}.  In each plot we
report the logarithm of the estimated stability time, $\log_{10} T$,
in color-code: blue (bottom of color legend) represents the smallest
stability time; yellow (top of color legend) refers to the largest
one.  We mark the actual position of Mercury, in the parameter space,
by a black dot. In addition, we draw a white curve through the
sub-space of parameters that lead to $\eps=2.06'$. This curve turns
out to be smooth and in perfect agreement with the one obtained
directly from~\eqref{myFor}. On its basis we also plot (dashed, black)
contour-curves corresponding to $\varepsilon=2.06'\pm5\%$ in order to
investigate the sensitivity of $\eps$ on the parameters.

The qualitative description of our results concerning the stability
time is as follows:
\begin{itemize}
\item the stability time is mostly influenced by the inclination as we
  can see from Figures~\ref{f:om2dot-incl}--\ref{f:om1dot-incl};

\item the stability time is only moderately influenced by changing the following 
parameters: the polar moment of
inertia (see Figures~\ref{f:crm-om2dot}--\ref{f:crm-om1dot}), the mean 
regression rate of the ascending node (see Figures~\ref{f:crm-om2dot} \& \ref{f:om1dot-om2dot}), 
and the eccentricity (see Figures~\ref{f:crm-ecc}--\ref{f:om1dot-ecc});

\item the precession rate of the perihelion argument does not seem to play 
a major role on the stability time as we can see in Figure~\ref{f:crm-om1dot} or 
\ref{f:om1dot-ecc}.

\end{itemize}

We now discuss in more detail the relations of the stability time,
$T$, and of the equilibrium obliquity, $\eps_*$, on specific system
parameters and only summarize in Table~\ref{tab:mercury_stabstat} the
outcome of our computations, where we report the ranges of the
effective stability time and obliquity, for all the different choices
of the parameters.

The stability times, in
Figures~\ref{f:om2dot-incl}--\ref{f:om1dot-incl}, increase from
$10^{12}$ to $10^{32}$ years for increasing inclination in the range
$0.05\leq i\leq 0.2$, and are only marginally influenced by the
parameters on the absciss{\ae}.  The obliquities turn out to increase for
increasing $|\dot\Omega|$ and $i$ in Figure~\ref{f:om2dot-incl} (the
largest value of $\eps$ is found in the upper left corner).  In
Figure~\ref{f:ecc-incl} the values for $\eps$ increase for decreasing
$e$ and increasing $i$ (being again largest in the upper left
corner). Obliquities do increase for increasing $c$ and $i$ in
Figure~\ref{f:crm-incl}, while $\eps$ stays constant for varying
$\dot\omega$ in Figure~\ref{f:om1dot-incl} but increases again for
increasing $i$.

The strong influence of the orbital inclination of Mercury is also
present in the obliquity ranges (see Table~\ref{tab:mercury_stabstat}):
for small $i$ we find $0.52'\leq\eps\leq 0.68'$ while for large $i$ we find 
$2.73'\leq\eps\leq4.72'$. The effect of the parameters
on the absciss{\ae} generally turn out to be small compared to changes in inclinations.
For the remaining cases the stability time ranges from $\log_{10} \ T=26.80$ (large
$c$ and small $e$), to $\log_{10} \ T=28.02$ (large $c$ and large $|\dot\Omega|$).
For these cases the smallest $\eps$ are still 2-3 times larger than for the
previous one (ranging from $1.27'$ to $1.77'$) while the maxima turn out to
be of the same order as before (ranging from $2.36'$ to $4.54'$).

We remark that the role of inclination $i$ on the stability time is
related to the role of $\eps$ and $K$ through the relationship
$\eps=K-i$. The application of formula \eqref{myFor} for constant
parameters of Mercury, but varying $i$, shows that with increasing $i$
we find increasing $\eps$. Making use of formula \eqref{myFor} again,
we find larger $\eps$ for larger values of $c$, $\dot\omega$ and
$\dot\Omega$ (keeping all remaining parameters constant). We carefully
checked in Figures~\ref{f:om2dot-incl}--\ref{f:om1dot-om2dot} that
maxima of the stability time correspond to maxima of obliquity $\eps$.
Moreover, we find that inclination $i$ most effectively increases
obliquity $\eps$ in \eqref{myFor} being consistent with our result.
Contrary, the role of $e$ on the stability time is different: we find
better stability for special values of $e$ within $0.085\leq e\leq
0.1$ (see Figures~\ref{f:crm-ecc}--\ref{f:om1dot-ecc}.  Moreover,
increasing $e$ gives smaller $\eps$ from \eqref{myFor}.  We conclude
that not only the role of $e$ on the realization of the 3:2 resonance
is special, but $e$ has a special role on the stability of the
resonance too.

Let us stress that the role of the eccentricity $e$ is quite subtle: a
non-zero value is needed in order to ensure the existence of a 3:2
resonance, but the eccentricity plays also the role of a perturbing
parameter.  In our study we assume that Mercury is placed in its
actual position, thus close to the 3:2 resonance, and just change the
value of the eccentricity, focusing only on the perturbation character
of the parameter.

\renewcommand{\arraystretch}{1.2}
\begin{table}
\begin{center}
\caption{Summary of stability times and obliquity ranges for different
parameter sets used to obtain Figure~\ref{f:om2dot-incl}--\ref{f:om1dot-ecc}}
\begin{tabular}{rcccc}
\hline
 set & $\varepsilon _{\min }[']$ & $\varepsilon _{\max }[']$ &
$\log_{10}T_{\min }[y]$ &   $\log_{10}T_{\max }[y]$ \\
\hline
 $\dot{\Omega }\times i$: & 0.52 & 4.72 & 12.00 & 32.00 \\
 $e\times i$: & 0.55 & 3.49 & 12.00 & 32.00 \\
 $c\times i$: & 0.59 & 3.13 & 12.00 & 32.00 \\
 $\dot{\omega }\times i$: & 0.68 & 2.73 & 12.00 & 32.00 \\
 $c\times\dot{\Omega }$: & 1.35 & 4.07 & 27.86 & 28.02 \\
 $c\times\dot{\omega }$: & 1.77 & 2.36 & 27.87 & 27.98 \\
 $\dot{\Omega }\times \dot{\omega }$: & 1.57 & 3.55 & 27.91 & 27.96 \\
 $c\times e$: & 1.42 & 3.02 & 26.80 & 28.00 \\
 $\dot{\Omega }\times e$: & 1.27 & 4.54 & 27.20 & 28.00 \\
 $\dot{\omega }\times e$: & 1.66 & 2.63 & 27.20 & 28.00 \\
\hline
\end{tabular}
\label{tab:mercury_stabstat}
\end{center}
\end{table}

\section{Conclusions and outlook}\label{sec:4}
We have investigated, analytically, the non-linear stability of Mercury's
3:2 spin-orbit resonance.  In particular, we have shown that Mercury is
currently placed in a very stable position in the parameter space.
Our study produces a strong bound on the longitudinal and
latitudinal librations over long, but finite times: we find that
libration widths up to $0.1$ rad stay bound for times exceeding the
age of the Universe.

However, from the results presented in the previous section, Mercury
does not seem to occupy the most stable configuration.  Indeed,
increasing the polar moment of inertia, $c$, and the inclination, $i$,
or increasing the mean regression rate of the ascending node
$|\dot\Omega|$ allows to get better estimates.

Our conclusions are valid on the basis of parameters that are relevant
for Mercury. The possible instability of spin-orbit resonances,
induced by non-linear perturbations, has been shown in, e.g.,
\cite{PavMac-2003}, \cite{BreEtal-2005}, and \cite{CelVoy-2010}. In
the latter study the authors find that the motion close to the 3:2
resonance, in the symmetric case ($A=B$), is essentially regular,
while a chaotic layer may appear increasing the asymmetry parameter
($S\equiv(B-A)/C$) and the eccentricity $e$.

In our computations we have, at most, $S\simeq10^{-4}$ for $c=0.3$.
Thus, our analytic estimates of the effective stability time is in
complete agreement with the results in \cite{CelVoy-2010}.  In fact,
our value of $S$ is one order of magnitude smaller compared to the
chaotic region numerically determined.

The regularity of the motion of Mercury has also been confirmed in
\cite{PavMac-2003}, where the authors find a soft transition from a
stable periodic motion to the unstable one in the 3:2 spin-orbit
resonance.  The authors also give an upper estimate for the ellipticity
($d_{cr}=3(B-A)/C$) of Mercury: $d_{cr}\simeq0.19626$.  Again, in our
application, the ellipticity lies well beyond $d_{cr}$.

It is interesting to notice that, in our study, the role of the
precession rate of the perihelion argument, $\dot\omega$, and the
eccentricity, $e$, on the 3:2 resonance stability is smaller compared
to the one of $\dot\Omega$ and $c$.  This result should be worthwhile
to be investigated further.

Our study is based on a simplified model where the orbital
ellipse of Mercury is kept constant, but precessing in the
node and perihelion. Additional perturbations that may act on semi-major axis,
orbital eccentricity and inclination may induce further perturbations
on the rotational motion of Mercury.  Furthermore, it would be
interesting to extend our study to be able to include internal
effects.

\section*{Acknowledgments}
We thank Ben\^oit Noyelles for interesting comments and discussions.
This work took benefit from the financial support of the contract
Prodex CR90253 from BELSPO.

\bibliographystyle{mn2e}
\bibliography{biblio}{}

\label{lastpage}

\end{document}